\documentclass[12 pt]{amsart}
\usepackage{amssymb,pdfsync}
\usepackage{graphicx}
\usepackage[mathscr]{eucal}
\usepackage{dsfont}
\usepackage[usenames]{color}
\definecolor{Purple}{rgb}{.7,0.08,0.6} 
 
\usepackage{enumerate}

\theoremstyle{plain}
\newtheorem{Thm}{Theorem}
\newtheorem{Cor}[Thm]{Corollary}
\newtheorem{Prop}[Thm]{Proposition}

\theoremstyle{definition}

\newcommand{\bndry}{b}
\newcommand{\B}{\mathbb B}

\newcommand{\C}{\mathbb C}
\newcommand{\Cin}{\mathbf C}
\newcommand{\Cino}{\mathbf C}
\newcommand{\Cint}{\Cin}
\newcommand{\Cinz}{\Cin_0}
\newcommand{\Cinmz}{\Cin_0}
\newcommand{\Cn}{\mathbb C^n}
\newcommand{\Ctro}{\mathcal C}

\newcommand{\Ctroe}{\mathbf C_\epsilon}

\renewcommand{\d}{\delta}
\newcommand{\dee}{\partial}
\newcommand{\deebar}{\overline\partial}

\newcommand{\De}{\Delta}
\newcommand{\Dmz}{D_0}
\newcommand{\Deo}{\De}
\newcommand{\Det}{\De}
\newcommand{\Deoe}{\De_\epsilon}
\newcommand{\Deme}{\De_\epsilon}
\newcommand{\Demz}{\De_0}
\newcommand{\Dez}{\Delta_0}
\newcommand{\Dme}{D_\epsilon}
\newcommand{\Do}{D}
\newcommand{\Dt}{D}
\newcommand{\Dz}{D_0}
\newcommand{\Dzm}{D_0}
\newcommand{\Dzme}{D_\epsilon}
\newcommand{\Dze}{D_\epsilon}

\newcommand{\fz}{f_0}
\newcommand{\fze}{f_\epsilon}

\newcommand{\la}{\lambda}
\newcommand{\lao}{\la}
\newcommand{\lat}{\la}
\newcommand{\laz}{\lambda_0}
\newcommand{\lamz}{\lambda_0}
\newcommand{\laoe}{\la_\epsilon}
\newcommand{\lame}{\la_\epsilon}
\newcommand{\Lo}{\Lambda}
\newcommand{\Loe}{\Lambda_\epsilon}

\newcommand{\mz}{m}

\renewcommand{\r}{\rho}
\newcommand{\ro}{\rho}
\newcommand{\roe}{\ro_\epsilon}
\newcommand{\rt}{\rho}
\newcommand{\rmo}{\rho}
\newcommand{\rz}{\rho_0}
\newcommand{\rmz}{\rho_0}
\newcommand{\rme}{\rho_\epsilon}
\newcommand{\R}{\mathbb R}

\newcommand{\si}{\sigma}

\newcommand{\sit}{\sigma}
\newcommand{\siz}{\sigma_0}
\newcommand{\simz}{\sigma_0}

\newcommand{\Smz}{\mathbf S_0}
\newcommand{\Smp}{\mathbf S_0'}
\newcommand{\Se}{\mathbf S_\epsilon}
\newcommand{\Sep}{\mathbf S_\epsilon'}
\newcommand{\Sz}{\mathbf S_0}
\newcommand{\Szp}{\mathbf S_0'}

\newcommand{\tze}{\tau_\epsilon}
\newcommand{\tme}{\tau_\epsilon}
\newcommand{\Tze}{\tau_\epsilon}

\begin{document}
\title
[Examples Cauchy-Leray]{The Cauchy-Leray integral:\\
counter-examples to the $L^p$-theory}
\date{\today}
\author[Lanzani and Stein]{Loredana Lanzani$^*$
and Elias M. Stein$^{**}$}
\thanks{$^*$ Supported in part by the National Science Foundation, award
DMS-1503612.}
\thanks{$^{**}$ Supported in part by the National Science Foundation, award
DMS-1265524.}
\address{
Dept. of Mathematics,       
Syracuse University 
Syracuse, NY 13244-1150 USA}
  \email{llanzani@syr.edu}
\address{
Dept. of Mathematics\\Princeton University 
\\Princeton, NJ   08544-100 USA }
\email{stein@math.princeton.edu}
\begin{abstract} We prove the optimality of the hypotheses guaranteeing the $L^p$-boundedness for the Cauchy-Leray integral in $\Cn$, $n\geq 2$, obtained in \cite{LS-4}.

Two domains, both elementary in nature, show that the geometric requirement of strong $\C$-linear convexity, together with regularity of order 2, are both necessary.
\end{abstract}
\maketitle

\section{Introduction}

The purpose of this paper is to present counter-examples that show that the $L^p$ results for the Cauchy-Leray integral obtained in \cite{LS-4} are essentially optimal. A second paper in this series will deal with counter-examples for the Cauchy-Szeg\H o projection, relevant to \cite{LS-5}.

Recall that in the case of the unit ball $\mathbb B$ in $\Cn$, the Cauchy-Leray integral and the 
Cauchy-Szeg\H o projection agree, and the same holds for the unbounded realization $U_0$ of $\mathbb B$. However when $n\geq 2$ for more general domains, these two operators are quite different. For the former, the result obtained in \cite{LS-4} states:
\medskip

{\em Suppose $D$ is a bounded domain in $\Cn$ whose boundary is of class $C^{1,1}$ and which is strongly $\C$-linearly convex. Then the induced Cauchy-Leray transform is bounded on $L^p(\bndry D, d\sigma)$ for $1<p<\infty$.}
\medskip

Our counter-examples show that these two restrictions on the boundary of $D$, the geometric condition of strong $\C$-linear convexity, and the regularity of degree 2, are in fact necessary.

We consider the following simple bounded domains in $\mathbb C^2$. Here $z_j=x_j+i\,y_j$, $j=1, 2$. The first domain is defined by
\begin{equation}\label{E:01}
|z_2|^2 + x_1^2 + y_1^4\ <\ 1\, .
\end{equation}

This is the domain given up to a translation by \eqref{E:6} in Section \ref{S:2}.  It has a
$C^\infty$ (in fact real-analytic) boundary, is strongly pseudo-convex, but not strongly $\C$-linearly convex.

The second example is
\begin{equation}\label{E:02}
|z_2|^2 + |x_1|^m + y_1^2\ <\ 1\, ,
\quad \text{where}\quad 1<m<2.
\end{equation}

This domain is given (up to a translation) by \eqref{E:24} in Section \ref{S:5}.  The domain has a boundary of class $C^{2-\epsilon}$, ($\epsilon = 2-m$), is strongly $\C$-linearly convex (and hence strongly pseudoconvex).

We show that for each domain the induced Cauchy-Leray transform is not bounded on $L^p$ for any $p$, $1\leq p\leq\infty$.

The idea of the proof is as follows: we assume that $L^p$-boundedness holds for one of these domains. Then an appropriate scaling and limiting argument shows that this positive result implies a corresponding conclusion in a limiting model domain, where it is much easier to supply  an explicit counter-example.

For the first example the limiting domain is the unbounded half-space $2\, \mathrm{Im}\, z_2\, >\, x_1^2$ (called
$\Dz$ in \eqref{E:1ab}), which is holomorphically equivalent to the more familiar half-space 
$U_0$, $2\, \mathrm{Im} z_2\, >\, |z_1|^2$, which itself is holomorphically equivalent to the unit ball. Note that the last two are strongly $\C$-linearly convex, but $\Dz$ is not. The hint that one might be led to a counter-example for $\Dz$, and then for the domain \eqref{E:01} is that its Cauchy-Leray operator is not ``pseudo-local''; 
(see \eqref{E:2} which shows that the kernel is singular away from the diagonal). 

The analysis of the second domain, represented by \eqref{E:02}, is parallel to that of the first domain.
For example, the corresponding limiting domain is $2\, \mathrm{Im}\, z_2\,>|x_1|^m$, see
 \eqref{E:25}.
 
A remark is in order about the family of boundary measures $\{d\mu_a\}_a$ given in \eqref{E:25a} that define
 the $L^p$ spaces, 
$L^p=L^p(\bndry D, d\mu_a)$ in the above results. Three examples of $d\mu_a$ are significant in many circumstances. First, the induced Lebesgue measure $d\sigma$; second, the Leray-Levi measure 
$d\lambda$, see \eqref{E:1a}; and third, the Fefferman measure
 \cite{B2}, \cite{F-1}, \cite{G}. 
  When the domain is smooth and strongly pseudo-convex the three measures give the same $L^p$ spaces, so that the counter-example holds for the domain \eqref{E:01} for all such measures.
In example \eqref{E:02} these measures are essentially different, yet the counter-example still
holds in all cases.

We should also call the reader's attention to the earlier relevant work in \cite{BaLa} where counter-examples are given for Cauchy-Leray integrals. However the less explicit and more complex nature 
of their construction and proof limit the results to the case $p=2$. It should be stressed that when 
$n\geq 2$ in general the Cauchy-Leray transform is far from ``self-adjoint'' and so failure of $L^2$-boundedness does not imply failure for any $p$, $p\neq 2$.

\section{The Cauchy-Leray integral; the model domain $\Dz$}\label{S:1}
For a bounded domain (say $C^2$) domain $D$ in $\Cn$, with a defining function 
$\r$, the corresponding Cauchy-Leray integral is
\begin{equation}\label{E:1}
\Cin (f) (z) =\int\limits_{\bndry D}\frac{1}{\De(w, z)^n}\, f(w)\, d\la (w)\, \quad z\in D.
\end{equation}
Here $f$ is (say) a bounded function on $\bndry D$, 
$$\De (w, z) = \langle\dee\r (w), w-z\rangle\, =\, \sum\limits_{j=1}^n\frac{\dee\r (w)}{\dee w_j} (w_j-z_j),$$ and
\begin{equation}\label{E:1a}
d\la (w) = \frac{j^*}{(2\pi i)^n}\left( \dee\r \wedge (\deebar\dee\r)^{n-1}\right )
\end{equation}
is the Leray-Levi measure (with $j^*$ the pull-back under the inclusion: $\bndry D\hookrightarrow \Cn$). We have $d\la (w) = \Gamma (w)\, d\si (w)\, ,$ with $d\si$ the induced Lebesgue
 measure on $\bndry D$ and
$$
\Gamma (w) = \frac{(n-1)!}{4\pi ^n} \frac{\mathcal L(w)}{|\nabla \r (w)|}\, ,
$$
where $\mathcal L(w)$ is the determinant of the Levi-form acting on the tangent space at $w\in\bndry D$. (See e.g., \cite[Ch. 7]{Ra}.) 

It is worth recalling two intrinsic properties of \eqref{E:1}:
\begin{itemize}
\item Its invariance under changes of coordinates that are given either by translations or unitary mappings of the space $\Cn$, see \cite{Bo}.
\item The independence of the Cauchy-Leray integral \eqref{E:1} of the particular defining function $\r$ of $D$.
\end{itemize}

We consider first the unbounded ``model domain'' $\Dz$ in $\mathbb C^2$.

\noindent With $z=(z_1, z_2)$, $z_j=x_j + i\, y_j$, it is defined by
\begin{equation}\label{E:1ab}
\Dz =\{z\, :\ 2\,\mathrm{Im} z_2 > x_1^2\}\, ,
\end{equation}
which is to be compared with the more familiar form of $\Dz$ given by 

\noindent $\{z\, :\ 2\,\mathrm{Im} z_2 > |z_1|^2\}$. These two domains are biholomorphically
equivalent via the mappings $z_1 \leftrightarrow z_1$, $z_2\leftrightarrow z_2 \pm i\, z_1^2$. Now the complex tangent space of these domains at the origin is the subspace $\{(z_1, 0)\}$. So since $|z_1|^2 = x_1^2 + y_1^2$ , and on $\mathbb R^2$ this is positive definite, this implies that the second domain is strongly $\C$-linearly convex. However because the form $x_1^2$ is degenerate along the direction $y_1$,
it follows that strong $\C$-linear convexity fails for the domain \eqref{E:1a}. (For more about these convexities see \cite{APS}; \cite{Ho}; \cite[Sect. 3.3]{LS-3}.)

With $\rz (z) = 2y_2-x_1^2$, and $\Dez (w, z) = \langle\dee\rz (w), w-z\rangle$, and if we write $w=(w_1, w_2)$, $w_j=u_j+i\, v_j$, $j=1, 2$, a simple calculation gives that
 with $w$ and $z\in \bndry \Dz$,
 
 \begin{equation}\label{E:2}
 \Dez(w, z) = \frac12 \big((u_1-x_1)^2 + 2i\, (u_1(v_1-y_1) + u_2 - x_2\big).
 \end{equation}

Now for fixed $z\in\bndry \Dz$, observe that $\Dez (w, z)$ vanishes on the one-dimensional variety given by $u_1=x_1$, $u_1(v_1-y_1) + u_2= x_2$. Also the Leray-Levi measure $d\laz$ corresponding to $\rz$ is $du_1\, dv_1\, du_2/(4\pi^2)\approx d\siz$,
 the induced Lebesgue measure  on $\bndry\Dz$, if we are near the origin.

To construct a counterexample for the model domain $\Dz$ we first choose a small constant $a$, which we will keep fixed throughout ($a=1/12$ will do). Then for any
$\d$,\ $0<\d<1$, we define:
\begin{equation}\label{E:3}
\left\{
\begin{array}{lll}
\mathbf S & = & \left\{|u_1|\leq a\d^2\, ,\ |v_1|\leq\frac12\, , \ |u_2|\leq a\d^2\right\} \\
\\
\mathbf S' & = & \left\{\d\leq |x_1|\leq 2\d\, ,\ |y_1|\leq\frac12\, , \ |x_2|\leq a\d^2\right\}
\end{array}
\right.
\end{equation}
with $w\in\bndry\Dz$ written as $(u_1+i\, v_1, u_2 +i\, u^2_1/2)$ and 
$z\in\bndry\Dz$ written as $(x_1+i\, y_1, x_2 +i\, x^2_1/2)$ and $\Sz$ and $\Szp$ 
being the corresponding sets for $w\in \bndry\Dz$, $z\in\bndry\Dz$ via \eqref{E:3}, then it follows that $\Dez (w, z)\neq 0$ whenever $w\in\Sz$ and $z\in\Szp$. Also
$$
\siz (\Sz)\approx \laz (\Sz) \approx \mz (\mathbf S)\approx \d^4\, \quad\mathrm{while}\quad
\siz (\Szp)\approx \laz (\Szp) \approx \mz (\mathbf S')\approx \d^3
$$
where $\siz$ denotes the induced Lebesgue measure on $\bndry\Dz$ and $\mz$ is the standard Euclidean measure on the three-dimensional Euclidean space.

We shall now test the presumed inequality
\begin{equation}\label{E:4}
\|\Cinz (f)\|_{L^p(\Szp, d\siz)}\leq A_p\, \|f\|_{L^p(\bndry\Dz, d\siz)}\, ,
\end{equation}
when $f$ is assumed to have support in $\Sz$. Under these circumstances
$\Cinz (f) (z)$ is well-defined as an absolutely convergent integral
$$
\int\limits_{\Sz}\frac{1}{|\Dez(w, z)|^2}\,f(w)\, d\laz (w)\, ,
$$
when $z\in\Szp$, in view of the non-vanishing of $\Dez (w, z)$ for these $w$ and $z$.
\begin{Prop}\label{P:1}
The inequality \eqref{E:4} fails for every $p$,\ $1\leq p\leq \infty$.
\end{Prop}
The constant $A_p$ is of course assumed to be independent of $f$. 

In fact it is clear by \eqref{E:2} and \eqref{E:3} that we have:
$$
\mathrm{Re}\, \Dez(w, z)\geq \frac{\d^2}{4}\, \quad \mathrm{and}\quad
|\mathrm{Im}\,\Dez(w, z)|\leq a\d^2 + 2a\d^2=3a\d^2.
$$

Then $\left(\mathrm{Re}\,\Dez (w, z)\right)^2 \geq 2\, |\mathrm{Im}\,\Dez(w, z)|^2$,
(which holds if $a\leq 1/12$). But $|\Dez(w, z)|\leq c\d^2$ for some $c>0$, thus
$$
\mathrm{Re}\frac{1}{\left(\Dez (w, z)\right)^2}\geq c'\d^4\, .
$$
Now take $f=\chi_{\Sz}$, the characteristic function of $\Sz$. Therefore 
$$\mathrm{Re}\big(\Cinz (f)\big)(z)\geq c'\d^{-4}\laz (\Sz)\geq c>0\,, \quad
 \mathrm{for }\quad z\in\Szp\, .$$
 So 
 $$
 \|\Cinz (f)\|^p_{L^p(\Szp)}\geq c\laz (\Szp) \geq c' \d^3\, ,
 $$
while 
$$
\|f\|^p_{L^p(\bndry\Dz)} \leq \laz(\Sz)\leq c\d^4\, .
$$
Since $\d^3$ is not $O(\d^4)$ for small $\d$, \eqref{E:4} cannot hold, 
when $p<\infty$. The case $p=\infty$ requires a separate but simpler argument 
which we omit.
\section{The first counter-example}\label{S:2}
We turn to the domain \eqref{E:01} in $\mathbb C^2$ and it is useful to consider a 
translate of it, given by
\begin{equation}\label{E:6}
\Do =\{|z_2-i|^2 + x_1^2 + y_1^4<1\}\, .
\end{equation}
We will show after rescaling and a passage to the limit, that we can reduce consideration of $\Do$ to $\Dz$. From \eqref{E:6} it is clear that 
$$
\ro (z) = x_1^2 + y_1^4 + x_2^2 + y_2^2 -2y_2
$$
is a defining function for $\Do$, and that $\ro$ is strongly pluri-subharmonic, 
so $\Do$ is strongly pseudo-convex. Moreover, since each of the four one-variable functions, $x_1^2$, $y_1^4$, $x_2^2$, $y_2^2 -2y_2$ are strictly convex, the domain $\Do$ is itself strictly convex. We note that 
$$
\mathrm{Re}\,\De (w, z) =
\mathrm{Re}\langle\dee\ro (w), w-z\rangle =
\frac12\big(\nabla\ro(w), w-z\big)_{\R}
$$
where $(\cdot, \cdot)_\R$ is the real inner product induced on $\R^4=\C^2$ from 
$\langle\cdot, \cdot\rangle$. With $w=(w_1, w_2)$, $w_j-u_j+i\, v_j\, ,$ $j=1, 2$, we claim that
\begin{equation}\label{E:6'}
\big(\nabla\ro(w), w-z\big)_{\R}\geq 
\end{equation}
$$
\geq (x_1-u_1)^2 + (x_2 -u_2)^2 + (y_2-v_2)^2 + (v_1^2+y_1^2)(v_1-y_1)^2,
$$
when $w, z\in\bndry\Do$.

To prove \eqref{E:6'} we use the identity
$$
f(\beta) -f(\alpha) = (\beta -\alpha)f'(\alpha) +
\int\limits_\alpha^\beta (\beta-\alpha) f''(t)\, dt
$$
for the functions $f_1=u_1^2$, $f_2=u_2^2$, $f_3= v_2^2 -2v_2$, and $f_4 = v_1^4$.
Similarly, we replace $w=(u_1+i\, v_1, u_2 +i\, v_2)$ by 
$z=(x_1+i\, y_1, x_2 +i\, y_2)$ and add the corresponding identities. Taking into account that $\ro(w) = f_1+ f_2+f_3+f_4$ wih $\ro(w)=0$ and the similar fact for 
$\ro (z)$, together with the observations that
$$
\int\limits_\alpha^\beta (\beta-\alpha) f''(t)\, dt = (\beta-\alpha)^2\quad \mathrm{if}\quad
f''(t)=2\, ,
$$
and
$$
\int\limits_\alpha^\beta (\beta-\alpha) f''(t)\, dt\geq (\beta^2+\alpha^2)(\beta-\alpha)^2,
\quad \mathrm{if}\quad f''(t)=12\, t^2\, ,
$$
then yields \eqref{E:6'}. 

Now \eqref{E:6'} shows that if $w$ and $z\in\bndry\Do$, $w\neq z$, then $z$ lies on one side of the (real) tangent  plane to $\bndry\Do$
 at $w$. By convexity of $\Do$, the same holds for $z\in \overline\Do \setminus \{w\}$.

Turning to the Cauchy-Leray integral of $\Do$ we recall two preliminary facts. First 
$\Cino (f)(z)$ is holomorphic in $z\in \Do$, if $f$ is say an integrable function on 
$\bndry\Do$. Second, whenever $F$ is a holomorphic function in $\Do$ which
is continuous in $\overline\Do$, and $f=F\big|_{\bndry\Do}$, then $\Cino(f)(z)=F(z),\ z\in\Do$. The latter fact follows from the Cauchy-Fantappi\`e formalism
 (see \cite{Ra} and \cite{LS-3}).

 Note that when $f$ is a bounded function supported in a closed set $\mathbf S$ in $\bndry D$, 
 then $\Cin (f)(z)$ is well-defined as a convergent integral whenever $z$ is outside the support
 $\mathbf S$. So certainly the extendability of $\Cin$ to a bounded operator on $L^p$ would imply 
 \begin{equation} \label{E:12}
 \|\Cin (f)\|_{L^p(\mathbf S')}\leq A_p\|f\|_{L^p}
 \end{equation}
 whenever $\mathbf S'$ is disjoint from $\mathbf S$, with $A_p$ independent of $f$, $\mathbf S$, and $\mathbf S'$.
 \begin{Thm}\label{T:1}
 For any $p$, $1\leq p\leq \infty$, the presumed inequality \eqref{E:12} fails when tested for a 
 bound $A_p$ independent of $f$ and its support $\mathbf S$, and with $\mathbf S'$ disjoint 
 from $\mathbf S$.
 \end{Thm}
   A further statement along these lines is made possible by the following fact (whose proof will
 appear elsewhere): whenever $f$ is of class $C^1$ on the boundary $\bndry D$, $\Cin (f)$
 is extendable to a continuous function on $\overline D$. With this we can define the Cauchy-Leray transform $\Ctro (f)$ for such $f$, by $\Ctro (f) = \Cin (f)\big|_{\bndry D}$.
 \begin{Cor} The mapping $f\mapsto \Ctro (f)$, initially defined for $C^1$ functions $f$ is not extendable to a bounded operator on $L^p(\bndry D)$.
 \end{Cor}
 
 \section{Proof of Theorem \ref{T:1}}\label{S:4}

 We shall obtain a contradiction to \eqref{E:12} by a scaling argument that leads us back to
 Proposition \ref{P:1}. 
 
 We define the maps $\tze$, $\epsilon >0$, on $\C^2$ by $\tze (z_1, z_2) = (\epsilon z_1, \epsilon^2 z_2)$ and set $\Dze =\tau_{\epsilon^{-1}}(\Do)$, with $\Do$ as in \eqref{E:6}. Then
 $\ro (z) = x_1^2 + y_1^4 + x_2^2 + y_2^2 -2y_2$ is a defining function for $\Do$ and 
 $\roe (z) = \epsilon^{-2}\ro (\tze (z))$ is a defining function for $\Dze$. Note that
 $$
 \roe (z) = x_1^2-2y_2 + \epsilon^2(y_1^4 + x_2^2 + y_2^2)\, ,
 $$
 from which it is clear that the domains $\Dze$ increase as $\epsilon$ decreases, with limit
 our model domain $\Dz =\{\rz (z)<0\}$ and $\rz (z) = x_1^2 -2y_2$.
 
 Observe also that
 $$
 \roe (z)\to \rz (z)\, ,\quad \mathrm{and}\quad \nabla\roe (z)\to \nabla \rz (z)
 $$
 uniformly on compact subsets of $\C^2$.
 
 We next define the ``transported'' measures $d\laoe$ on $\bndry\Dze$ by
 \begin{equation}\label{E:16}
 \int\limits_{z\in\bndry\Dze}\!\!\! F(z)\, d\laoe (z) \ =\ 
 \epsilon^{-4}\!\!\! \int\limits_{z\in\bndry\Dz}\!\!\! F(\tau_{\epsilon^{-1}}(z))\, d\lao (z)
 \end{equation}
 for any continuous functions $F$ defined on $\bndry\Dze$.
 
  Here $d\lao$ is the Leray-Levi measure on $\bndry\Do$, but note that $d\laoe$ is not the Leray-Levi measure on $\bndry\Dze$: as a result the operator $\Ctroe$ defined below is
  {\em not} the Cauchy-Leray integral of $\bndry\Dze$. We can also define the corresponding action of $\tze$ on functions by
  $$
  \Tze (F)(z) = F(\tze (z)) = F\circ\tze\, .
  $$
 Then by what has been said above, $\Tze$ maps  $L^p(\bndry\Do, d\lao)$ to 
 $L^p(\bndry\Dze, d\laoe)$ and we have the ``isometry'':
 \begin{equation}\label{E:17}
 \|F\|_{L^P(\bndry\Do, d\lao)} \ =\ \epsilon^{4/p}\, \|\,\Tze (F)\|_{L^p(\bndry\Dze, d\laoe)}\, .
 \end{equation}
  Now let $\Deoe (w, z) = \langle\dee\roe (w), w-z\rangle$. Then 
  $$\Deo (\tze (w), \tze(z)) = 
  \langle\dee\ro (\tze w), \tze(w-z)\rangle =
  \epsilon^2 \langle \dee\roe (w), w-z\rangle = 
  \epsilon^2 \Deoe (w, z)\, .
  $$
  Hence
  \begin{equation}\label{E:18}
  \Deoe (w, z)\ =\ \epsilon^{-2}\,\Deo (\tze (w), \tze (z)).
  \end{equation}
  We now define the operator $\Ctroe$, by setting
  \begin{equation}\label{E:19}
  \Ctroe (f) (z) \ =\int\limits_{\bndry\Dze}\!\! \frac{1}{\Deoe (w, z)^2}\, f(w)\, d\laoe (w)
  \end{equation}
  which is well-defined as the integral above for any bounded function $f$, as long as $z$ is outside the support of $f$. Our next claim is that 
  \begin{equation}\label{E:20}
  \Ctroe (f) (z) \ =\ \Ctro (f\circ\tau_{\epsilon^{-1}})(\tze (z)),
  \end{equation}
  for bounded $f$ on $\bndry\Dze$, if $z\in\bndry\Dze$ lies outside the support of $f$. In fact going back to the definition \eqref{E:1}, and using \eqref{E:16} and \eqref{E:18}, we see that
  $$
  \Cin (f\circ \tau_{\epsilon^{-1}})(\tze (z)) \ =\int\limits_{\bndry\Do}
  \!\! \frac{1}{\Deo (w, \tze (z))^2}\, f(\tau_{\epsilon^{-1}}(w))\, d\lao (w)\ =
  $$
 $$
 \epsilon^4 \int\limits_{\bndry\Dze}
  \!\! \frac{1}{\Deo (\tze (w), \tze (z))^2}\, f(w)\, d\laoe (w)\ =
  \int\limits_{\bndry\Dze}
  \!\! \frac{1}{\Deoe (w, z)^2}\, f(w)\, d\laoe (w)\, 
 $$ 
 showing \eqref{E:20}.
 
 Next, if \eqref{E:12} held, then by \eqref{E:17} we would also have
 \begin{equation}\label{E:21}
 \|\Ctroe (F)\|_{L^p(\Sep, d\laoe)}\leq A_p \|F\|_{L^p(\bndry\Dze, d\laoe)}
 \end{equation}
 whenever $F$ is a bounded function on $\bndry\Dze$ and $\Sep$ is a closed subset of
 $\bndry\Dze$, disjoint from the support of $F$. 
 
 At this point we restrict the attention to the unit ball $\B$ in $\C^2$, and we exploit the common coordinate system for $\bndry\Dz \cap \B$, and $\bndry\Dze\cap \B$, when $\epsilon$ is small.
 That is in this ball, $\bndry\Dz$ is the graph over $(x_1+i\,y_1, x_2)$ given by
 $(x_1+i\,y_1, x_2 + i\, x_1^2/2)$, while $\bndry\Dze$ is the graph given by

 \noindent $(x_1+i\,y_1, x_2 + i\, \left(x_1^2/2 + \Phi_\epsilon (x_1, y_1, x_2)\right))$, with
  $\Phi_\epsilon (x_1, y_1, x_2) = O(\epsilon^2)$.
  
  Now if in these coordinates we write $d\lao = \Lo (x_1, y_1, x_2)\, dx_1\, dy_1\, dx_2$ and
  $d\lambda_\epsilon = \Loe (x_1, y_1, x_2)\, dx_1\, dy_1\, dx_2$, then by \eqref{E:16} we have
  \begin{equation}\label{E:22}
  \Loe (x_1, y_1, x_2) \ =\ \Lo (\epsilon\, x_1, \epsilon\,y_1, \epsilon^2 x_2)\, .
  \end{equation}
  Finally we take $f$ to be the function $\chi_{\mathbf S}$, the characteristic function of the set 
  $\mathbf S$ in \eqref{E:3}, and lift $f$ to functions on $\bndry\Dze$ and $\bndry\Dz$ respectively, by setting $\fze (z_\epsilon) = f(x_1, y_1, x_2)$ when
  $z_\epsilon = (x_1+i\,y_1, x_2 + i\, \left(x_1^2/2 + \Phi_\epsilon (x_1, y_1, x_2)\right))$, and 
  $\fz (z_0)= f(x_1+i\,y_1, x_2 + i\, x_1^2/2)$. 
  
  We also lift the sets $\mathbf S$ and $\mathbf S'$ to $\Se$ and $\Sep$ (subsets of $\bndry\Dze$) in the same way. Our claim is with that notation
  \begin{equation}\label{E:23}
  \Ctroe (\fze) (z_\epsilon)\ \to\ \Cinz (\fz) (z_0)\, ,\quad \mathrm{if}\ z_0\in\Szp\, .
  \end{equation}
  In fact, 
  $$\Ctroe (\fze) (z_\epsilon)\ =
  \int\limits_{\mathbf S}\frac{1}{\Deoe(w, z_\epsilon)^2}\, f(u_1, v_1, u_2)\, 
  \Loe (u_1, v_1, u_2)\, du_1\, dv_1\, du_2\, .
  $$
  However by \eqref{E:22} $\Loe (u_1, v_1, u_2)\to \Lo(0, 0, 0)$, and moreover
  $\Deoe(w, z_\epsilon)\to\Dez (w, z_0)$ because $\nabla\roe\to\nabla\rz$, while
   $\Dez (w, z_0)\neq 0$ if $w\in\Sz$ and $z_0\in\Szp$. This gives \eqref{E:23}.
   
   As a result \eqref{E:21} with $F=\fze$ leads to
   $$
   \|\Cinz (\fz)\|_{L^p(\Szp, d\laz)}\leq A_p \|\fz\|_{L^p(\bndry\Dz, d\laz)}
   $$
   which contradicts Proposition \ref{P:1}, proving the theorem.
   \section{The second counter-example}\label{S:5}
   Here the domain will be taken to be
   \begin{equation}\label{E:24}
   \Dt =\{|z_2-i|^2+|x_1|^m+y_2^2<1\}\, ,\quad \mathrm{with}\quad 1<m<2\, .
   \end{equation}
   Its model domain is
   \begin{equation}\label{E:25}
   \Dzm \ =\ \{2\,\mathrm{Im}\, z_2>|x_1|^m\}.
   \end{equation}
   For any $f$ that is bounded on $\bndry\Dt$, the Cauchy-Leray integral $\Cint (f) (z)$ is well-defined for $z$ that lies in the complement of the support of $f$. As in the previous sections we will show that the mapping $f\mapsto \Cint (f)$ fails to be bounded in $L^p$ in the sense that
   the proposed inequality 
   \begin{equation}\label{E:26}
   \|\Cint (f)\|_{L^p(\mathbf S', d\si)}\leq A_p\|f\|_{L^p(\bndry\Dt, d\sit)}
   \end{equation}
   cannot hold. Here $\mathbf S'$ is any set disjoint from the support of $f$, and the bound $A_p$ is assumed independent of $f$ and $\mathbf S'$.
   
   The proof of this assertion follows the same lines as in Sections \ref{S:1}-\ref{S:4} for the domain \eqref{E:6}, and so we will only discuss the minor differences that occur.
   
   The defining function of $\Dt$ is $\rt (z) = |x_1|^m + y_1^2+ x_2^2 + y_2^2 -2y_2$, and that 
   of $\Dmz$ is $\rmz (z) = |x_1|^m -2y_2$. Note that both domains are of class $C^{2-\alpha}$,
   with $\alpha = 2-m$. Also since $|x_1|^m\,,\ y_1^2\,, \ x_2^2\, ,\ y_2^2-2y_2$ are strongly 
   convex functions of one variable, the domain $\Dt$ is strongly convex, hence strongly $\C$-linearly convex, and in fact strongly pseudoconvex in the following sense: the domain $D$ is
   exhausted by an increasing family of smooth domains $\{D_\gamma\}_\gamma$ which are uniformly strongly pseudo-convex, with defining function:
    $\rho_\gamma (z) = (x_1^2+\gamma)^{m/2} + y_1^2 + x_2^2 + y_2^2 -2y_2$.

   This convexity implies that $\mathrm{Re}\, \Det (w, z)>0$ for $w\in\bndry\Dt$ and $z\in\overline\Dt$, except when $z=w$.
   
   Returning to the model domain, if $\Demz (w, z) = \langle\dee\rmz (w), w-z\rangle$ a calculation gives
   \begin{equation}\label{E:27}
   \Demz (w, z) = 
   \end{equation}
   $$
 =  |x_1|^m - |u_1|^m + m\,[u_1]^{m-1}(u_1-x_1) \ +
   \ i\, \big(m\,[u_1]^{m-1}(v_1-y_1) + 2(x_2-y_2)\big)\, .
   $$
   Here we have used the notation 
   $$
 [u_1]^{m-1}\ =\ \frac1m\,  \frac{d}{du_1}|u_1|^m \ =\ |u_1|^{m-1}\,\mathrm{sign}\, u_1\, .
   $$
   Now we set
   \begin{equation}\label{E:28}
\left\{
\begin{array}{lll}
\mathbf S & = & \left\{|u_1|\leq a\,\d^2\, ,\ |v_1|\leq\d^{2-m}\, , \ |u_2|\leq \d^m\right\} \\
\\
\mathbf S' & = & \left\{\d\leq |x_1|\leq 2\d\, ,\ |y_1|\leq\d^{2-m}\, , \ |x_2|\leq \d^m\right\}
\end{array}
\right.
\end{equation}
and let $\Smz$ and $\Smp$ be the corresponding induced sets on $\bndry\Dmz$.

We have that near the origin 
$$
d\lamz \approx |u_1|^{m-2}\, du_1\, dv_1\, du_2,\quad \mathrm{and}\quad
d\simz \approx du_1\, dv_1\, du_2\, ,
$$
where $d\lamz$ and $d\simz$ are the Leray-Levi measure and the induced Lebesgue measure
on $\bndry\Dmz$. Thus
\begin{equation}\label{E:29}
\lamz (\Smz)\ \approx \ \d^{2m}
\end{equation}
\begin{equation}\label{E:30}
\simz (\Smz)\ \approx\ \d^4\, ,\quad\mathrm{and}\quad \simz(\Smp)\ \approx\ \d^3\, .
\end{equation}
By \eqref{E:27} it follows that 
$$
\mathrm{Re}\,\Demz (w, z) \ \gtrsim \ \d^m\, ,\quad\mathrm{while}\quad |\mathrm{Im}\,\Demz (w, z)|\lesssim \d^m\, ,
$$
and if we choose $a$ sufficiently small, then
$$
\mathrm{Re}\,\frac{1}{\big(\Demz (w, z)\big)^2}\geq c\, \d^{-2m}\quad \mathrm{if}\ w\in\Smz\, ,\quad \mathrm{while}\ z\in\Smp\, .
$$
Taking $\fz$ to be the characteristic function of $\Smz$, we therefore get 
$\mathrm{Re}\, \Cinmz (\fz)(z)\geq c>0$, for $z\in\Smp$. Hence this gives a contradiction to \eqref{E:26} in the case when $\Cinmz$ is the Cauchy-Leray integral of the model domain $\Dmz$.

To pass to the domain $\Dt$ we carry out the scaling via
\begin{equation}\label{E:31}
\tme (z_1, z_2) \ =\ (\epsilon z_1, \epsilon^m z_2)
\end{equation}

The domain $\Dzme = \tau_{\epsilon^{-1}}(\Dt)$ has a defining function
$$\rme (z)\ =\,  \epsilon^{-m}\rmo (\tme (z))\ =\ \rmz (z) + \epsilon^m(x_2^2+y_2^2) + \epsilon^{2-m}y_1^2\, ,$$
which converges to $\rmz (z) = |x_1|^m - 2y_2$. 

We also set
$\Deme (w, z) = \langle \dee\rme (w), w-z\rangle$. The transported measure $d\lame$ on 
$\Dme$ is now defined by the identity
\begin{equation}\label{E:32}
\int\limits_{\bndry\Dme} F(z)\, d\lame (z) \ =\ \epsilon^{-2m}\int\limits_{\bndry\Dt} F(\tau_{\epsilon^{-1}}(z))\, d\lat (z)\, ,
\end{equation}
(compare with \eqref{E:16}).

 Finally, if we assumed that \eqref{E:26} held for the domain $\Dt$ we can then see by the reasoning in Section \ref{S:4} that the corresponding result would hold for the model domain
 $\Dmz$ achieving our desired contradiction.
 
 We point out in closing that \eqref{E:26} fails not only for the $L^p$ norms taken with the induced Lebesgue measure but others as well. 
 Consider the measures $d\mu_a$ on $\bndry\Dt$
 given by
 \begin{equation}\label{E:25a}
 d\mu_a =\left(\frac{\mathcal L}{|\nabla\rho|}\right)^{\!\! a} \! d\sit\, ,
 \end{equation}
 with $\mathcal L$ the determinant of the Levi-form, and $d\sit$ the induced Lebesgue measure.
 
 Then the argument above
  shows that \eqref{E:26} fails when
 $$
-\infty< a <\frac{1}{2-m}\, ,
 $$
(because $\mathcal L^a\approx |x_1|^{(m-2)a}$).
 Here the factor $\epsilon^{-2m}$ in \eqref{E:32} is replaced with $\epsilon^{a(2-m)-2-m}$ to reflect the new transported measures $d\mu_{a,\epsilon}$.
 
 The case $a=0$ corresponds to induced Lebesgue measure; the case $a=1$, to the Leray-Levi measure; and the case $a=1/3$ to the Fefferman measure \cite{B2}, \cite{F-1}, \cite{G}.

\end{document}